\documentclass[12pt]{article}
\usepackage{latexsym}
    \title{{\bf  Genus-zero modular functors and 
intertwining operator algebras}}
    \author{Yi-Zhi Huang}
    \date{}
    \begin{document}
    \bibliographystyle{alpha}
    \maketitle

    \input amssym.def
    \input amssym
    \newtheorem{rema}{Remark}[section]
    \newtheorem{propo}[rema]{Proposition}
    \newtheorem{theo}[rema]{Theorem}
   \newtheorem{defi}[rema]{Definition}
    \newtheorem{lemma}[rema]{Lemma}
    \newtheorem{corol}[rema]{Corollary}
     \newtheorem{exam}[rema]{Example}
\newcommand{\binom}[2]{{{#1}\choose {#2}}}
	\newcommand{\nno}{\nonumber}
	\newcommand{\lbar}{\bigg\vert}
\newcommand{\mbar}{\mbox{\large $\vert$}}
	\newcommand{\p}{\partial}
	\newcommand{\dps}{\displaystyle}
	\newcommand{\bra}{\langle}
	\newcommand{\ket}{\rangle}
 \newcommand{\res}{\mbox{\rm Res}}
\renewcommand{\hom}{\mbox{\rm Hom}}
 \newcommand{\pf}{{\it Proof.}\hspace{2ex}}
 \newcommand{\epf}{\hspace{2em}$\Box$}
 \newcommand{\epfv}{\hspace{1em}$\Box$\vspace{1em}}
\newcommand{\nord}{\mbox{\scriptsize ${\circ\atop\circ}$}}
\newcommand{\wt}{\mbox{\rm wt}\ }
\newcommand{\clr}{\mbox{\rm clr}\ }

\begin{abstract} 
In \cite{H2} and \cite{H4}, the author introduced the notion of
intertwining operator algebra, a nonmeromorphic generalization of the notion
of vertex operator algebra involving monodromies. The problem of constructing
intertwining operator algebras from representations of
suitable vertex operator algebras were solved implicitly 
earlier in  \cite{H1}. In the
present paper, we generalize the geometric and operadic formulation of
the notion of vertex operator algebra given in \cite{H-1}, \cite{H0},
\cite{HL1}, \cite{HL2} and \cite{H3}
to the notion of intertwining
operator algebra. We show that the category of intertwining operator
algebras of central charge $c\in {\Bbb C}$ is isomorphic to the
category of  algebras over rational genus-zero
modular functors (certain analytic partial operads) of central charge
$c$ satisfying certain generalized meromorphicity. 
This result is one main step in the 
construction of genus-zero conformal field theories from representations of 
vertex operator algebras announced in \cite{H2}.
One byproduct of the proof of the present isomorphism
theorem is a geometric construction of (framed) braid group
representations from intertwining operator algebras and thus from
representations of suitable vertex operator algebras.
\end{abstract}

\renewcommand{\theequation}{\thesection.\arabic{equation}}
\renewcommand{\therema}{\thesection.\arabic{rema}}
\setcounter{equation}{0}
\setcounter{rema}{0}
\setcounter{section}{-1}

\section{Introduction}

In this paper, we show that  the category of
intertwining operator algebras of central charge $c\in {\Bbb C}$ is
isomorphic to the category of 
algebras over
rational genus-zero modular functors of central charge $c$ satisfying a 
certain generalized meromorphicity. This result is one main step in the 
construction of genus-zero conformal field theories from representations of 
vertex operator algebras announced in \cite{H2}.

The geometric and operadic
formulation of the notion of vertex operator algebra (\cite{H-1}, \cite{H0},
\cite{HL1}, \cite{HL2}, \cite{H3})
establishes a direct connection between vertex operator algebras and
the geometry of genus-zero Riemann surfaces with extra structures. 
In particular, 
this formulation gave a simple and conceptual definition of vertex 
operator algebra: A vertex operator algebra is a meromorphic algebra
over a vertex partial operad. The equivalence
between this formulation and the algebraic formulation of the
notion of vertex operator algebra (\cite{B}, \cite{FLM}) provides a practical
way to construct genus-zero holomorphic conformal field theories. See
\cite{H3} for a presentation of this theory.

On the other hand, vertex operator algebras are not enough for either
the study of conformal field theories or mathematical problems
related to vertex operator algebras.  Most conformal field theories
have monodromies which are described by multivalued operator-valued
functions. In the study of the moonshine module vertex operator
algebra constructed Frenkel, Lepowsky and Meurman \cite{FLM}
and related problems, it is often necessary to study intertwining
operators. The study of these multivalued operator-valued functions and
intertwining operators leads us naturally to the notion of
intertwining operator algebra introduced in \cite{H2} and \cite{H4}. A
construction of these algebras from representations of
suitable vertex operator algebras were given implicitly in
\cite{H1}. In fact, even for a problem whose statement involes only
vertex operator algebras, the solution often involves intertwining
operator algebras. One example is the conceptual construction of 
the vertex operator algebra structure on the moonshine module given in 
\cite{H1.5}.

In this paper, we generalize the geometric and operadic formulation of
the notion of vertex operator algebra to a geometric and operadic
formulation of the notion of intertwining operator algebra.  We
construct a rational genus-zero modular functor (a certain analytic
partial operad) from an intertwining operator algebra and show that
the intertwining operator algebra gives an algebra satisfying a
certain generalized meromorphicity over this partial operad.  Since a
generalized-meromorphic algebra over a rational genus-zero modular
functor is a genus-zero weakly holomorphic conformal field theory
(see \cite{S1}, \cite{S2} and \cite{H2}), we obtain
a construction of genus-zero weakly holomorphic conformal field
theories from intertwining operator algebras. In particular, we obtain
a geometric construction of (framed) braid group representations from
vertex operator algebras and modules.  Conversely, we also show that
an algebra satisfying the generalized meromorphicity over such a
partial operad (a generalized-meromorphic algebra over a rational
genus-zero modular functor) gives an intertwining operator
algebra. The main result of the present paper is that the category of
intertwining operator algebras of central charge $c\in {\Bbb C}$ is
isomorphic to the category of canonically-generalized-meromorphic
algebras over
rational genus-zero modular functors of central charge $c$ (a
subcategory of the category of generalized-meromorphic
algebras over
rational genus-zero modular functors of central charge $c$).  In
particular, we obtain a simple and conceptual definition of
intertwining operator algebra: An intertwining operator algebra is a
(canonically-)generalized-meromorphic algebra 
over a rational genus-zero modular
functor.

The method used in the present paper is almost the same as in
\cite{H3}, except that the duality properties for vertex operator
algebras are replaced by the more complicated duality properties for
intertwining operators. In fact, the constructions and the proofs in
the present paper are adaptions or modifications of those in
\cite{H3}.  Thus, in many steps, instead of giving detailed
constructions or proofs, we shall indicate the differences between the
discussions in the present paper and those in \cite{H3}, and refer the
reader to the corresponding constructions or proofs in 
\cite{H3} for more details.

We assume that the reader is familiar with the material in \cite{H3}. 
In Section 1, we give  definitions of genus-zero modular
functor and rational genus-zero modular
functor. 
We also construct flat
connections on holomorphic vector bundles underlying rational genus-zero
modular functors and show that rational genus-zero modular functors give
representations of the (framed) braid groups
in this section. A definition of intertwining operator algebra is
given in Section 2. In Section 3, we construct a rational genus-zero
modular functor and a generalized-meromorphic algebra over it from an
intertwining operator algebra.
The main result of the present paper is given in Section 4.

\paragraph{Acknowledgments} 
I would like to express my gratitude to Professors Gu Chaohao,
Hu Hesheng, Hong Jiaxing and Dr. Zhou Zixiang 
for their hospitality. I am grateful to Jim Lepowsky for comments.
This research is supported in part by NSF
grant  DMS-9622961.

\renewcommand{\theequation}{\thesection.\arabic{equation}}
\renewcommand{\therema}{\thesection.\arabic{rema}}
\setcounter{equation}{0}
\setcounter{rema}{0}
\section{Genus-zero modular functors}

We first give definitions of genus-zero modular functor and rational
 genus-zero modular functor in terms of the language of partial
 operads. Modular functors were first introduced by G. Segal in his
 geometric formulation of the notion of conformal field theory (see
 \cite{S1} and \cite{S2}). See also \cite{H2} for the genus-zero case.
The definitions given in this section are
slightly different from the one in \cite{S1} and \cite{S2} (restricted to
the genus-zero case) because we are using them to prove
precise theorems. 
We formulate the definition of modular
 functor in the genus-zero case more generally
 to include the irrational case, and more precisely for our study of
intertwining operator algebras. 
We have added Axioms 3 and 5 in Definitions \ref{mf} and 
 and \ref{rmf}, respectively, and we have also modified Axioms 2 and 6.  
These modifications are necessary for formulating the precise 
geometric definition
of intertwining operator algebra and our main theorem.

\begin{defi}\label{mf}
{\rm A {\it genus-zero modular functor} is
an analytic partial
operad ${\cal M}$ together with a morphism $\pi:
{\cal M}\to K$ of analytic partial
operads a finite-rank holomorphic vector bundle structure 
on the triple 
$({\cal M}(j), K(j), \pi)$ for any $j\in {\Bbb N}$
satisfying the following axioms:

\begin{enumerate}

\item For any ${\cal Q}\in {\cal M}(k)$, ${\cal Q}_{1}\in
{\cal M}(j_{1}), \dots, {\cal Q}_{k}\in {\cal M}(j_{k})$, $k,
j_{1}, \dots, j_{k}\in {\Bbb N}$, the substitution
$\gamma_{{\cal M}}({\cal Q}; {\cal Q}_{1},\dots, {\cal Q}_{k})$ in
${\cal M}$ exists if (and only if) 
$$\gamma_{K}(\pi({\cal Q});
\pi({\cal Q}_{1}), \dots, \pi({\cal Q}_{k}))$$ 
exists.\label{a1}

\item Let $Q\in K(k)$, $Q_{1}\in K(j_{1}), \dots, Q_{k}\in K(j_{k})$,
$k, j_{1}, \dots, j_{k}\in {\Bbb N}$, such that $\gamma(Q; Q_{1},
\dots, Q_{k})$ exists.  The map {}from the Cartesian product of the
fibers over $Q$, $Q_{1}, \dots, Q_{k}$ to the fiber over
$\gamma_{K}(Q; Q_{1}, \dots, Q_{k})$ induced {}from the substitution map
of ${\cal M}$ is multilinear and surjective.

\item The rank of ${\cal M}(0)$ is $1$.

\end{enumerate}}
\end{defi}

The simplest examples of genus-zero modular functors are 
${\Bbb C}$-extensions of $K$ discussed in Section 6.9 of \cite{H3}. 
In fact any genus-zero modular functor ${\cal L}$ such that ${\cal L}(j)$
for any $j\in {\Bbb N}$ is a line bundle is a ${\Bbb C}$-extension of $K$.

Given a  genus-zero modular functor ${\cal M}$, 
let ${\cal M}^{*}(j)$ be the dual vector bundle of ${\cal M}(j)$ for
$j\in {\Bbb N}$. Then it is clear that 
${\cal M}^{*}=\{{\cal M}^{*}(j)\}_{j\in {\Bbb N}}$
 is a genus-zero modular functor. We call ${\cal M}^{*}$ the {\it dual}
of ${\cal M}$.  Given genus-zero modular functors ${\cal M}_{1}$ and
${\cal M}_{2}$, Let $({\cal M}_{1}\otimes {\cal M}_{2})(j)$ be the
tensor product bundle ${\cal M}_{1}(j)\otimes {\cal M}_{2}(j)$ of ${\cal
M}_{1}(j)$ and ${\cal M}_{2}(j)$.  Then it is also clear that ${\cal
M}_{1}\otimes {\cal M}_{2}=\{({\cal M}_{1}\otimes {\cal
M}_{2})(j)\}_{j\in {\Bbb N}}$ is a genus-zero modular functor. 
We call ${\cal
M}_{1}\otimes {\cal M}_{2}$ the {\it tensor product} of ${\cal M}_{1}$ 
and ${\cal M}_{2}$.

\begin{defi}\label{rmf}
{\rm A {\it rational genus-zero modular functor}
is a genus-zero modular functor ${\cal M}$ and a
finite set ${\cal A}$ satisfying the following axioms:

\begin{enumerate}

\setcounter{enumi}{3}

\item For any $j\in {\Bbb N}$ and $\alpha_{0}, \alpha_{1}, \dots, 
\alpha_{j}\in {\cal 
A}$, there are finite-rank holomorphic vector bundles ${\cal 
M}_{\alpha_{1}\cdots \alpha_{j}}^{\alpha_{0}}(j)$ over $K(j)$ such that 
$${\cal 
M}(j)=\oplus_{\alpha_{0}, \alpha_{1}, \dots, \alpha_{j}\in {\cal A}}{\cal 
M}_{\alpha_{1}\cdots \alpha_{j}}^{\alpha_{0}}(j)$$ 
where $\oplus$ denotes the direct
sum operation for vector bundles.

\item There exists $\epsilon\in {\cal A}$ such that 
the rank of ${\cal M}^{\alpha}(0)$ is $1$ if $\alpha=\epsilon$
and is $0$ if $\alpha_{0}\ne \epsilon$. For any $\alpha_{0}, 
\alpha_{1}\in {\cal A}$, the rank of ${\cal M}^{\alpha_{0}}_{\alpha_{1}}(1)$
is $1$ if  $\alpha_{0}=\alpha_{1}$ and is $0$ if $\alpha_{0}\ne \alpha_{1}$.

\item \label{last} For any $j\in {\Bbb N}$ and $\beta_{0}, \beta_{1}, 
\dots, \beta_{j}\in {\cal A}$,
let $P^{\beta_{0}}(j)$ be the projection {}from ${\cal M}(j)$ to
$\oplus_{\alpha_{1}, \dots, \alpha_{j}\in {\cal A}}
{\cal  M}_{\alpha_{1}\cdots
\alpha_{j}}^{\beta_{0}}(j)$ and $P_{\beta_{1}\cdots \beta_{j}}(j)$ the 
projection {}from
${\cal M}(j)$ to $\oplus_{\alpha_{0}\in {\cal A}} {\cal  M}_{\beta_{1}\cdots
\beta_{j}}^{\alpha_{0}}(j)$. Then for any $k\in {\Bbb Z}_{+}$, $j_{1}, \dots,
j_{k}\in {\Bbb N}$,
$$\gamma_{{\cal M}}
= \sum_{\beta_{1}, \dots, \beta_{k}\in {\cal
A}}\gamma_{{\cal M}}\circ (P_{\beta_{1}\cdots \beta_{k}}(k)
\times P^{\beta_{1}}(j_{1})\times \cdots P^{\beta_{k}}(j_{k})),$$
and for any
$Q\in K(k)$, $Q_{1}\in K(j_{1}), \dots, Q_{k}\in K(j_{k})$,
$k\in {\Bbb Z}_{+}$, 
$j_{1}, \dots, j_{k}\in {\Bbb N}$, such that $\gamma(Q; Q_{1},
\dots, Q_{k})$ exists and any $\beta_{0}, \alpha_{l}^{(i)}\in
{\cal A}$, $l=1, \dots, j_{i}$, $i=1, \dots, k$,  
the linear map {}from 
$$\oplus_{\beta_{1}, \dots, \beta_{k}\in {\cal A}}
{\cal M}^{\beta_{0}}_{\beta_{1}, \dots
\beta_{k}}\mbar_{Q}\otimes {\cal M}^{\beta_{1}}_{\alpha^{(1)}_{1} \cdots
\alpha^{(1)}_{j_{1}}}\mbar_{Q_{1}}\otimes \cdots
\otimes 
{\cal M}^{\beta_{k}}_{\alpha^{(k)}_{1} \cdots
\alpha^{(k)}_{j_{k}}}\mbar_{Q_{k}}$$
(where we use ${\cal M}\mbar_{Q}$ to denote the fiber of a vector 
bundle ${\cal M}$ at an element $Q$ of the base manifold)
to ${\cal M}\mbar_{\gamma_{K}(Q; Q_{1}, \dots, Q_{k})}$ 
induced from the map
$\gamma_{{\cal M}}\circ (P_{\beta_{1}\cdots \beta_{k}}(k)
\times P^{\beta_{1}}(j_{1})\times \cdots P^{\beta_{k}}(j_{k}))$ 
is an isomorphism.
\end{enumerate}}
\end{defi}

\begin{rema}
{\rm Recall the notions of rescaling group of a partial operad and
rescalable partial operad in \cite{HL1}, \cite{HL2} and Appendix C of
\cite{H3}. In \cite{H3}, it was shown that ${\Bbb C}^{\times}$ can be
identified with a rescaling group of $K(1)$ and that $K$ is a ${\Bbb
C}^{\times}$-rescalable partial operad.  From Axiom \ref{a1} and \ref{last}
in the
definitions above, we see that the set 
${\cal M}^{\times}(1)|_{{\Bbb C}^{\times}}$
of nonzero elements of 
the fibers over elements of the rescaling group
of $K$ is a rescaling group of ${\cal M}$ and
${\cal M}$ is an ${\cal M}^{\times}(1)|_{{\Bbb C}^{\times}}$-rescalable 
partial
operad.}
\end{rema}

We now construct
a canonical holomorphic flat connection $\nabla(j)$ on ${\cal M}(j)$ for each
$j\in {\Bbb N}$. Recall the partial suboperad $\widehat{K}$
and $\overline{K}$ discussed
in \cite{H3}: For any $j\in {\Bbb N}$,
$\widehat{K}(j)$ is the subset  of $K(j)$ consisting of elements of the form
$$(z_{1}, \dots, z_{j-1}; A(z_{j}; 1), (a_{0}^{(1)}, {\bf 0}),
\dots, (a_{0}^{(j)}, {\bf 0})),$$
where $z_{1}, \dots, z_{j-1}\in {\Bbb C}^{\times}$ satisfying
$z_{k}\ne z_{l}$ for $k\ne l$, $z_{j}\in {\Bbb C}$,
$a_{0}^{(1)}, \dots,
a_{0}^{(j)}\in {\Bbb C}^{\times}$,
and $A(z_{j}, 1)$ is the sequence of complex numbers 
whose first component is $1$ and whose other components are $0$.
Then 
$$\widehat{K}=\{\widehat{K}(j)\}_{j\in {\Bbb N}}
$$
is a partial suboperad of
$K$.
For any $j\in {\Bbb N}$, $\overline{K}(j)$
is the subset of $\widehat{K}(j)$ consisting of elements of the form
$$(z_{1}, \dots, z_{j-1}; A(z_{j}; 1), (1, {\bf 0}),
\dots, (1, {\bf 0})).$$
Then $\overline{K}=\{\overline{K}(j)\}_{j\in {\Bbb N}}$
 is also a partial suboperad of $K$. Note that for any
$j\in {\Bbb Z}_{+}$, $\overline{K}(j)$
can be identified with the configuration space
$$\{(z_{1}, \dots, z_{j})\in {\Bbb C}^{j}\; |\; z_{k} \ne z_{l}, k\ne
l\}$$ 
by the map
$$(z_{1}, \dots, z_{j-1}; A(z_{j}; 1), (1, {\bf 0}),
\dots, (1, {\bf 0}))\mapsto (z_{1}-z_{j}, \dots, z_{j-1}-z_{j}, -z_{j}).$$

We consider the restrictions $\widehat{{\cal M}}(j)$
of the holomorphic vector bundles ${\cal M}(j)$ to $\widehat{K}(j)$
for $j\in {\Bbb N}$. Then $\widehat{{\cal M}}
=\{\widehat{{\cal M}}(j)\}_{j\in {\Bbb N}}$ is a partial suboperad of 
${\cal M}$. We first construct holomorphic flat connections on
 $\widehat{{\cal M}}(j)$ for $j\in {\Bbb N}$.

Let $U$ be a simply connected open subset of ${\Bbb C}^{\times}$ containing
$1$. We first construct for any $\alpha\in {\cal A}$ 
a holomorphic flat connection 
on the restriction ${\cal M}^{\alpha}_{\alpha}(1)
\mbar_{U}$ of ${\cal M}^{\alpha}_{\alpha}(1)$ to
$U\subset {\Bbb C}^{\times}\subset \widehat{{\cal M}}(1)$.
The composition map for the partial operad ${\cal M}$ gives a holomorphic
bundle map from ${\cal M}^{\alpha}_{\alpha}(1)
\mbar_{U}\times {\cal M}^{\alpha}_{\alpha}(1)
\mbar_{U}$ to ${\cal M}^{\alpha}_{\alpha}(1)
\mbar_{U\cdot U}$ where $U\cdot U$ is the image of $U\times U$ under the 
multiplication of complex numbers. Since $U$ is simply connected, we
can always find a holomorphic trivialization of 
${\cal M}^{\alpha}_{\alpha}(1)
\mbar_{U}$ such that any pair of pullbacks of the
flat sections of the product bundle over $U$ under the trivialization
is mapped to a section of ${\cal M}^{\alpha}_{\alpha}(1)
\mbar_{U\cdot U}$ whose restriction to $U$ is the pullback of 
a flat connection of the product bundle over $U$. This trivialization 
gives a holomorphic flat connection 
on ${\cal M}^{\alpha}_{\alpha}(1)
\mbar_{U}$. These holomorphic flat connections together give a holomorphic 
connection on ${\cal M}(1)\mbar_{U}$.

We now construct the flat connection on ${\cal M}(j)$ 
for any $j\in {\Bbb N}$.
Let 
$$Q=(z_{1}, \dots, z_{j-1}; A(z_{j}, 1), (a_{0}^{(1)}, {\bf 0}), 
\dots, (a_{0}^{(j)}, {\bf 0}))$$ 
be an element of $K(j)$ and 
\begin{eqnarray*}
\lefteqn{\ell=\{Q(t)=(z_{1}(t), \dots, z_{j-1}(t); 
A(z_{j}(t); 1), (a_{0}^{(1)}(t), {\bf 0}), \dots, 
(a_{0}^{(j)}(t), {\bf 0}))}\nno\\
&&\hspace{10em}\in \widehat{K}(j)\;|\;t\in [0, 1]\}\hspace{10em}
\end{eqnarray*}
be a smooth path beginning at $Q(0)=Q$. Then for sufficiently small $t$,
$(a_{0}^{(i)})^{-1}a_{0}^{(i)}(t)\in U$, $i=1, \dots, j$, and
\begin{eqnarray*}
Q(t)&=&((\cdots (Q_{^{1}}\infty_{^{0}}(A(z_{1}(t)-z_{1}; 1), 
((a_{0}^{(1)})^{-1}a_{0}^{(1)}(t), {\bf 0})))
_{^{2}}\infty_{^{0}}\cdots )_{^{j-1}}\infty_{^{0}}\nno\\
&&\hspace{1em}
(A(z_{j-1}(t)-z_{j-1}; 1), ((a_{0}^{(j-1)})^{-1}a_{0}^{(j-1)}(t), {\bf 0}))
)_{^{j}}\infty_{^{0}}\nno\\
&&\hspace{2em}
(A(z_{j}(t); 1), ((a_{0}^{(j)})^{-1}a_{0}^{(j)}(t), {\bf 0})).
\end{eqnarray*}

From the definition of genus-zero modular functor, the substitution
maps for ${\cal M}$ give an isomorphism from the tensor product of the
fibers of $\widehat{\cal M}(j)$, $\widehat{\cal M}(1), \dots,
\widehat{\cal M}(1)$, $\widehat{\cal M}(1)$ over $Q$,
$(A(z_{1}(t)-z_{1}; 1), (a_{0}^{(1)}(t), {\bf 0})), \dots,
(A(z_{j-1}(t)-z_{j-1}; 1), (a_{0}^{(j-1)}(t), {\bf 0})),$
$(A(z_{j}(t); 1), (a_{0}^{(j-1)}(t), {\bf 0}))$, respectively, to the
fiber of $\widehat{\cal M}(j)$ over $Q(t)$. On the other hand, the
holomorphic 
flat connection ${\cal M}(1)\mbar_{U}$ gives isomorphisms from the fibers of
${\cal M}(1)\mbar_{U}$ over $$(A(z_{i}(t)-z_{i}; 1), 
((a_{0}^{(i)})^{-1}a_{0}^{(i)}(t), {\bf 0})), i=1, \dots, j-1,$$  
and $$(A(z_{j}(t); 1), (a_{0}^{(j)}(t), {\bf 0}))$$ to the fiber
of ${\cal M}(1)\mbar_{U}$ over $({\bf 0}, (1, {\bf 0}))$. Thus we
obtain an isomorphism from the fiber of $\widehat{\cal M}(j)$ over $Q$
to the fiber of $\widehat{\cal M}(j)$ over $Q(t)$.  Repeating the
procedure above starting from points in the path other than $Q$, we
obtain for any $t\in [0, 1]$ such an isomorphism from the fiber of
$\widehat{\cal M}(j)$ over $Q$ to the fiber of $\widehat{\cal M}(j)$
over $Q(t)$.  Given any element of the fiber of $\widehat{\cal M}(j)$
over $Q$, these isomorphisms give a path in $\widehat{\cal M}(j)$
beginning at this element such
that its projection image in $\widehat{K}(j)$ is the given path
$\ell$. Since $Q$ and $\ell$ are arbitrary, we obtain a connection  on
$\widehat{\cal M}(j)$. From the construction, we see that this
connection is flat and holomorphic.  

Since $K(j)$ is homotopically
equivalent to $\widehat{K}(j)$, we can extend the holomorphic vector
bundle $\widehat{\cal M}(j)$ over $\widehat{K}(j)$ and
the flat holomorphic connection 
just constructed trivially to a holomorphic vector
bundle ${\cal M}^{\sim}(j)$ over $K(j)$ and a holomorphic flat
connection on ${\cal M}^{\sim}(j)$, respectively. We can also extend
the partial operad structure on $\widehat{\cal M}$ to the sequence
${\cal M}^{\sim}= \{{\cal M}^{\sim}(j)\}_{j\in {\Bbb N}}$ so that
${\cal M}^{\sim}$ becomes a genus-zero modular functor.

Consider the genus-zero modular functor ${\cal M}\otimes 
({\cal M}^{\sim})^{*}$.
Since the fibers of the restriction of $({\cal M}\otimes 
({\cal M}^{\sim})^{*})(j)$ to $\widehat{K}(j)$ 
for any $j\in {\Bbb N}$ are tensor products of vector spaces and 
their dual spaces, we have a canonical isomorphism from this restriction
to a trivial line bundle over $\widehat{K}(j)$. Since 
$K(j)$ is homotopically equivalent to $\widehat{K}(j)$, this canonical
isomorphism can be extended to a canonical isomorphism from 
$({\cal M}\otimes 
({\cal M}^{\sim})^{*})(j)$ to a trivial line bundle
${\cal L}(j)$ over $K(j)$. 
The genus-zero modular functor structure on 
$({\cal M}\otimes 
({\cal M}^{\sim})^{*})(j)$ induces a genus-zero modular functor structure
on ${\cal L}=\{{\cal L}(j)\}_{j\in {\Bbb N}}$ and the canonical isomorphisms
constructed above give a morphism of genus-zero modular functors. 
So ${\cal M}$ is isomorphic to ${\cal M}^{\sim}\otimes {\cal L}$.
Since ${\cal L}(j)$, $j\in {\Bbb N}$, are all line bundles, ${\cal L}$ is
a ${\Bbb C}$-extension of $K$. By Theorem D.6.3 in \cite{H3},
there exists a complex number $c$ such that ${\cal L}$ is isomorphic to 
the $c$-th power $\tilde{K}^{c}$ of $\tilde{K}$. Thus ${\cal M}$ is isomorphic 
to ${\cal M}^{\sim}\otimes \tilde{K}^{c}$. Since for 
any $j\in {\Bbb N}$, there are canonical
holomorphic flat connections on ${\cal M}^{\sim}$ and 
$\tilde{K}^{c}(j)$, we obtain a holomorphic flat connection
$\nabla(j)$ on 
${\cal M}(j)$. 

In the construction above, we obtain a complex number $c$. We call this 
complex number the {\it central charge} of the genus-zero modular functor 
${\cal M}$.

We now consider the restrictions  $\widehat{\cal M}$ and
$\overline{\cal M}$ of a genus-zero modular functor 
${\cal M}$  to
the partial suboperad $\widehat{K}$ and 
$\overline{K}$, respectively. Recall that for any $j\in {\Bbb N}$
$\overline{K}(j)$ can identified with the configuration space
$F(j)$. So we can view $\overline{\cal M}(j)$ as a vector bundle over
$F(j)$ with a holomorphic flat connection. Since by the definition
of genus-zero modular functor, the action of 
$S_{j}$ on $\overline{\cal M}$ 
are given by isomorphisms of holomorphic vector 
bundles covering the action of $S_{j}$ on $F(j)$, 
the holomorphic vector bundle $\overline{\cal M}(j)$ on $F(j)$ and
its holomorphic flat connection induce a holomorphic vector bundle
on $F(j)/S_{j}$ and a holomorphic flat connection on this induced bundle.
In fact this holomorphic vector bundle and the holomorphic flat connection 
on it is the restriction of a holomorphic vector bundle
on $\widehat{\cal M}/S_{j}$ and a holomorphic flat connection on this 
vector bundle induced from 
the holomorphic vector bundle $\overline{\cal M}(j)$ on $F(j)$ and
its holomorphic flat connection.
We know that a flat connection on a vector bundle over a connected
manifold gives a structure of a representation of the fundamental
group of the base manifold on the fiber of the vector bundle over any point
on the manifold.  We also know that the braid group $B_{j}$ on $n$
strings is by definition the fundamental group of $F(j)/S_{j}$ and
the framed braid group on $j$ strings 
is the fundamental group $\widehat{\cal M}/S_{j}$. 
So we obtain:

\begin{theo}\label{braid}
Let ${\cal M}$ be a genus-zero modular functor, $j\in {\Bbb Z}_{+}$ and
${\cal V}$ the fiber over any point in $F(j)\subset K(j)$ of the vector
bundle ${\cal M}(j)$. Then ${\cal V}$ has a natural structure of a
representation of the framed braid group  on $j$ strings. In particular, 
${\cal V}$ has a natural structure of a
representation of the braid group $B_{j}$ on $j$ strings.
\end{theo}

\renewcommand{\theequation}{\thesection.\arabic{equation}}
\renewcommand{\therema}{\thesection.\arabic{rema}}
\setcounter{equation}{0}
\setcounter{rema}{0}
\section{Intertwining operator algebras}

The following definition of intertwining operator algebra 
was first given in \cite{H2}:

\begin{defi}
{\rm 
An {\it intertwining operator algebra of central charge
$c\in {\Bbb C}$} consists of the following data:

\begin{enumerate}

\item A finite-dimensional commutative associative algebra $A$ and a
basis ${\cal A}$ of $A$ containing the identity $\epsilon\in A$ such that all
the structure constants ${\cal  N}_{\alpha_{1}\alpha_{2}}^{\alpha_{3}}$, 
$\alpha_{1},
\alpha_{2}, \alpha_{3}\in {\cal A}$, are in ${\Bbb N}$.

\item  A vector space
$$W=\coprod_{\alpha\in
{\cal A}, n\in {\Bbb R}}W^{\alpha}_{(n)},\; \mbox{\rm for}\; 
w\in W^{\alpha}_{(n)},\;
n=\wt w,\; \alpha=\clr w$$
doubly graded
by ${\Bbb R}$ and ${\cal A}$
 (graded by {\it weight} and by {\it color}, respectively).

\item For each triple $(\alpha_{1}, \alpha_{2}, \alpha_{3})\in {\cal A}\times
{\cal A}\times {\cal A}$, an ${\cal
N}_{\alpha_{1}\alpha_{2}}^{\alpha_{3}}$-dimensional 
subspace
${\cal V}_{\alpha_{1}\alpha_{2}}^{\alpha_{3}}$ 
of the vector space of all linear maps
$W^{\alpha_{1}}\otimes W^{\alpha_{2}}\to W^{\alpha_{3}}\{x\}$, or equivalently,
 an ${\cal N}_{\alpha_{1}\alpha_{2}}^{\alpha_{3}}$-dimensional vector space
${\cal V}_{\alpha_{1}\alpha_{2}}^{\alpha_{3}}$ whose elements are linear maps
\begin{eqnarray*}
{\cal Y}: W^{\alpha_{1}}&\to& \mbox{Hom}(W^{\alpha_{2}}, W^{\alpha_{3}})\{x\}\\
w&\mapsto &{\cal Y}(w, x)
=\sum_{n\in {\Bbb R}}{\cal Y}_{n}(w)x^{-n-1} \;(\mbox{where}\;
{\cal Y}_{n}(w)\in \mbox{End}\;W).
\end{eqnarray*}

\item Two distinguished vectors ${\bf 1}\in W^{\epsilon}$ ({\it the vacuum})
and $\omega\in W^{\epsilon}$ ({\it the Virasoro element}). 

\end{enumerate}

\noindent These data satisfy the
following axioms for $\alpha_{1}, \alpha_{2}, 
\alpha_{3}, \alpha_{4}, \alpha_{5}, 
\alpha_{6}\in {\cal A}$, 
$w_{(\alpha_{i})}\in W^{\alpha_{i}}$, $i=1, 2, 3$, and $w_{(\alpha_{4})}'
\in W'_{\alpha_{4}}$:
\begin{enumerate}

\item The {\it weight-grading-restriction conditions}: 
For any $n\in {\Bbb Z}$ and $\alpha\in {\cal A}$, 
$$\dim W^{\alpha}_{(n)}< \infty$$
and 
$$W^{\alpha}_{(n)}=0$$
for $n$
sufficiently small.

\item Axioms for intertwining operators:

\begin{enumerate}

\item The {\it single-valuedness condition}: for any ${\cal Y}\in {\cal 
V}_{\epsilon\alpha_{1}}^{\alpha_{1}}$, 
$${\cal Y}(w_{(\alpha_{1})}, x)\in \mbox{Hom}(W^{\alpha_{1}},
W^{\alpha_{1}})[[x, x^{-1}]].$$

\item The {\it lower-truncation property for vertex operators}: for any
${\cal Y}\in {\cal V}_{\alpha_{1}\alpha_{2}}^{\alpha_{3}}$,
${\cal Y}_{n}(w_{(\alpha_{1})})w_{(\alpha_{2})}=0$ for $n$ sufficiently large.

\end{enumerate}

\item Axioms for the vacuum:

\begin{enumerate}

\item The {\it identity property}: for any ${\cal Y}\in {\cal 
V}_{\epsilon \alpha_{1}}^{\alpha_{1}}$, there is $\lambda_{{\cal Y}}
\in {\Bbb C}$ such that 
${\cal Y}({\bf 1}, x)=\lambda_{{\cal Y}}I_{W^{\alpha_{1}}}$, where
$I_{W^{\alpha_{1}}}$ on the right is the identity
operator on $W^{\alpha_{1}}$.

\item The {\it creation property}: for any ${\cal Y}\in {\cal 
V}_{\alpha_{1}\epsilon}^{\alpha_{1}}$, 
there is $\xi_{{\cal Y}}\in {\Bbb C}$ such that ${\cal 
Y}(w_{(\alpha_{1})}, x){\bf 1}\in W[[x]]$ and $\lim_{x\to 0}
{\cal Y}(w_{(\alpha_{1})},
x){\bf 1}=\xi_{{\cal Y}}w_{(\alpha_{1})}$ (that is, 
${\cal Y}(w_{(\alpha_{1})},
x){\bf 1}$ involves only nonnegative
integral powers of $x$ and the constant term is $\xi_{{\cal Y}}
w_{(\alpha_{1})}$).

\end{enumerate}

\item Axioms for products and iterates of intertwining operators:

\begin{enumerate}

\item The {\it convergence properties}: for any $m\in {\Bbb Z}_{+}$,
$\alpha_{i}, \beta_{i}, \mu_{i}\in {\cal A}$, $w_{(\alpha_{i})}
\in W^{\alpha_{i}}$, ${\cal Y}_{i}\in {\cal 
V}_{\alpha_{i}\;\beta_{i+1}}^{\mu_{i}}$, $i=1, \dots, m$, $w_{(\mu_{1})}'
\in (W^{\mu_{1}})'$ and 
$w_{(\beta_{m})}\in W^{\beta_{m}}$, the series
$$\langle w_{(\mu_{1})}', {\cal Y}_{1}(w_{(\alpha_{1})}, x_{1})
\cdots{\cal Y}_{m}(w_{(\alpha_{m})},
x_{m})w_{(\beta_{m})}\rangle_{W^{\mu_{1}}}\mbar_{x^{n}_{i}=e^{n\log z_{i}},
i=1, \dots, m, n\in {\Bbb R}}$$ 
is absolutely convergent when $|z_{1}|>\cdots >|z_{m}|>0$, and for 
any ${\cal Y}_{1}\in {\cal 
V}_{\alpha_{1}\alpha_{2}}^{\alpha_{5}}$ and ${\cal Y}_{2}\in {\cal 
V}_{\alpha_{5}\alpha_{3}}^{\alpha_{4}}$, the series
$$\langle w_{(\alpha_{4})}', {\cal Y}_{2}({\cal Y}_{1}(w_{(\alpha_{1})}, 
x_{0})w_{(\alpha_{2})},
x_{2})w_{(\alpha_{3})}\rangle_{W^{\alpha_{4}}}
\mbar_{x^{n}_{0}=e^{n\log (z_{1}-z_{2})},
x^{n}_{2}=e^{n\log z_{2}}, n\in {\Bbb R}}$$ is absolutely convergent when
$|z_{2}|>|z_{1}-z_{2}|>0$.

\item The {\it associativity}: for any ${\cal Y}_{1}\in {\cal 
V}_{\alpha_{1}\alpha_{5}}^{\alpha_{4}}$ and ${\cal 
V}_{\alpha_{2}\alpha_{3}}^{\alpha_{5}}$, there exist ${\cal Y}^{\alpha}_{3}
\in {\cal 
V}_{\alpha_{1}\alpha_{2}}^{\alpha}$ and ${\cal Y}^{\alpha}_{4}\in {\cal 
V}_{\alpha\alpha_{3}}^{\alpha_{4}}$ for all $\alpha\in {\cal A}$
such that the (multivalued) analytic function 
$$\langle w_{(\alpha_{4})}', 
{\cal Y}_{1}(w_{(\alpha_{1})}, x_{1}){\cal Y}_{2}(w_{(\alpha_{2})}, 
x_{2})w_{(\alpha_{3})}\rangle_{W}\mbar_{x_{1}=z_{1},
x_{2}=z_{2}}$$
on $\{(z_{1}, z_{2})\in {\Bbb C}\times {\Bbb C}\;|\; |z_{1}|>|z_{2}|>0\}$ 
and the (multivalued) analytic function
$$\sum_{\alpha\in {\cal A}}\langle w_{(\alpha_{4})}', {\cal Y}^{\alpha}_{4}
({\cal Y}^{\alpha}_{3}(w_{(\alpha_{1})},
x_{0})w_{(\alpha_{2})}, x_{2})w_{(\alpha_{3})}\rangle_{W^{\alpha_{4}}}
\mbar_{x_{0}=z_{1}-z_{2},
x_{2}=z_{2}}$$ 
on
$\{(z_{1}, z_{2})\in {\Bbb C}\times {\Bbb C}\;|\;
|z_{2}|>|z_{1}-z_{2}|>0\}$ are equal on $\{(z_{1}, z_{2})\in
{\Bbb C}\times {\Bbb C}\;|\;|z_{1}|> |z_{2}|>|z_{1}-z_{2}|>0\}$.

\end{enumerate}

\item Axioms for the Virasoro element:

\begin{enumerate}

\item The {\it Virasoro algebra relations}: Let $Y$ be the
element of ${\cal V}_{\epsilon\alpha_{1}}^{\alpha_{1}}$ such that $Y({\bf 1},
x)=I_{W^{\alpha_{1}}}$ and let $Y(\omega, x)=\sum_{n\in
{\Bbb Z}}L(n)x^{-n-2}$. Then
$$[L(m), L(n)]=(m-n)L(m+n)+\frac{c}{12}(m^{3}-m)\delta_{m+n, 0}$$
for $m, n\in {\Bbb Z}$.

\item The {\it $L(0)$-grading property}: 
$L(0)w=nw=(\wt w)w$ for $n\in {\Bbb R}$
and $w\in W_{(n)}$.

\item The {\it $L(-1)$-derivative property}: For any ${\cal Y}\in {\cal 
V}_{\alpha_{1}\alpha_{2}}^{\alpha_{3}}$, 
$$\frac{d}{dx}{\cal Y}(w_{(\alpha_{1})}, x)
={\cal Y}(L(-1)w_{(\alpha_{1})}, x).$$

\item The {\it skew-symmetry}: There is a linear map $\Omega$ {}from
${\cal V}_{\alpha_{1}\alpha_{2}}^{\alpha_{3}}$ 
to ${\cal V}_{\alpha_{2}\alpha_{1}}^{\alpha_{3}}$
such that for any ${\cal Y}\in {\cal V}_{\alpha_{1}\alpha_{2}}^{\alpha_{3}}$,
$${\cal Y}(w_{(\alpha_{1})}, x)w_{(\alpha_{2})}
=e^{xL(-1)}(\Omega({\cal Y}))(w_{(\alpha_{2})},
y)w_{(\alpha_{1})}\mbar_{y^{n}=e^{in\pi}x^{n}}.$$

\end{enumerate}
\end{enumerate}}
\end{defi}

To make our study slightly easier, we shall assume in the present
paper that intertwining operator algebras also satisfy the following
additional condition: For any $\alpha\in {\Bbb A}$, there exists
$h_{\alpha}\in {\Bbb R}$ such that $(W^{\alpha})_{(n)}=0$ for $n\not
\in h_{\alpha}+{\Bbb Z}$. But this is in fact a very minor 
restriction and in addition, all the results on intertwining operator
algebras satisfying this additional property can be generalized easily
to intertwining operator algebras not satisfying this condition.

The intertwining operator algebra defined above is denoted by 
$$(W,
{\cal A},  \{{\cal 
V}_{\alpha_{1}\alpha_{2}}^{\alpha_{3}}\}, {\bf 1}, \omega)$$ 
or simply by $W$. The
commutative associative algebra $A$ is called the {\it Verlinde
algebra}\index{Verlinde algebra} or the {\it fusion algebra of
$W$}\index{fusion algebra}. The linear maps in 
${\cal 
V}_{\alpha_{1}\alpha_{2}}^{\alpha_{3}}$ are called 
{\it intertwining operators of type 
${\alpha_{3}\choose \alpha_{1}\alpha_{2}}$}.

\begin{rema}
{\rm There are also other equivalent definitions of intertwining operator 
algebras using only formal variables. See \cite{H4} and \cite{H5} 
for details.}
\end{rema}

We shall also need the generalized rationality and commutativity 
for intertwining operator
algebras. Here we state a theorem which gives these properties
together with the associativity:

\begin{theo}\label{dual}
Let 
$$(W,
{\cal A},  \{{\cal 
V}_{\alpha_{1}\alpha_{2}}^{\alpha_{3}}\}, {\bf 1}, \omega)$$
be an intertwining operator algebra.
Then for any $\alpha_{1}, \alpha_{2}, \alpha_{3}, \alpha_{4},
\alpha_{5}\in {\cal A}$, any $w_{(\alpha_{1})}\in W_{\alpha_{1}}$, 
$w_{(\alpha_{2})}\in W^{\alpha_{2}}$, $w_{(\alpha_{3})}\in 
W^{\alpha_{3}}$, $w'_{(\alpha_{4})}\in (W^{\alpha_{4}})'$,
and any intertwining operators ${\cal Y}_{1}$ and ${\cal Y}_{2}$ of
type ${W^{\alpha_{4}}\choose W^{\alpha_{1}W\alpha_{5}}}$ and 
${W^{\alpha_{5}}\choose W^{\alpha_{2}W\alpha_{3}}}$, respectively,
there exist real numbers $r_{j}, s_{j}, t_{j}$, $j=1, \dots, p$, 
polynomial functions 
$g_{j}(z_{1}, z_{2})$, $j=1, \dots, p$, in $z_{1}$ and $z_{2}$, 
intertwining operators ${\cal Y}_{3}^{\alpha}$, ${\cal
Y}_{4}^{\alpha}$,  ${\cal Y}_{5}^{\alpha}$, ${\cal
Y}_{6}^{\alpha}$, $\alpha\in {\cal A}$, of type 
${W^{\alpha}\choose W^{\alpha_{1}W\alpha_{2}}}$,
${W^{\alpha_{4}}\choose W^{\alpha W\alpha_{3}}}$,
${W^{\alpha_{4}}\choose W^{\alpha_{2}W\alpha}}$,
${W^{\alpha_{5}}\choose W^{\alpha_{1}W\alpha_{3}}}$, respectively,
 such that 
$$\langle w_{(\alpha_{4})}', 
{\cal Y}_{1}(w_{(\alpha_{1})}, x_{1}){\cal Y}_{2}(w_{(\alpha_{2})}, 
x_{2})w_{(\alpha_{3})}\rangle_{W}\mbar_{x_{1}=z_{1},
x_{2}=z_{2}},$$
$$\sum_{\alpha\in {\cal A}}\langle w_{(\alpha_{4})}', {\cal Y}^{\alpha}_{4}
({\cal Y}^{\alpha}_{3}(w_{(\alpha_{1})},
x_{0})w_{(\alpha_{2})}, x_{2})w_{(\alpha_{3})}\rangle_{W^{\alpha_{4}}}
\mbar_{x_{0}=z_{1}-z_{2},
x_{2}=z_{2}}$$
and 
$$\sum_{\alpha\in {\cal A}}\langle w_{(\alpha_{4})}', 
{\cal Y}^{\alpha}_{5}(w_{(\alpha_{1})}, x_{1})
{\cal Y}^{\alpha}_{2}(w_{(\alpha_{2})}, 
x_{2})w_{(\alpha_{3})}\rangle_{W}\mbar_{x_{1}=z_{1},
x_{2}=z_{2}}$$
are the restrictions to the regions $|z_{1}|>|z_{2}|>0$,
$|z_{2}|>|z_{1}-z_{2}|>0$ and $|z_{2}|>|z_{1}|>0$, respectively, 
of the (multivalued) analytic
functions 
$$\sum_{j=1}^{p}
z_{1}^{-r_{j}}z_{2}^{-s_{j}}
(z_{1}-z_{2})^{-t_{j}}g_{j}(z_{1}, z_{2}).$$
\end{theo}

For a proof of this result, see Section 2 of \cite{H4}.

The theorem above can be generalized easily to more than two intertwining
operators. Here we state the generalized rationality part of this
generalization:

\begin{theo}\label{dual2}
Let 
$$(W,
{\cal A},  \{{\cal 
V}_{\alpha_{1}\alpha_{2}}^{\alpha_{3}}\}, {\bf 1}, \omega)$$
be an intertwining operator algebra. 
Then for any $m\in {\Bbb Z}_{+}$,
$\alpha_{i}, \beta_{i}, \mu_{i}\in {\cal A}$, $w_{(\alpha_{i})}
\in W^{\alpha_{i}}$, ${\cal Y}_{i}\in {\cal 
V}_{\alpha_{i}\;\beta_{i+1}}^{\mu_{i}}$, $i=1, \dots, m$, $w_{(\mu_{1})}'
\in (W^{\mu_{1}})'$ and 
$w_{(\beta_{m})}\in W^{\beta_{m}}$,
there exist real numbers $r^{(k)}_{j}, s^{(kl)}_{j}$, 
$j=1, \dots, p$, $k, l=1, \dots, m$, $k>l$, and
polynomial functions 
$g_{j}(z_{1}, , \dots, z_{m})$, $i=1, \dots, m$, in $z_{1}$ and $z_{2}$,
such that
$$\langle w_{(\mu_{1})}', {\cal Y}_{1}(w_{(\alpha_{1})}, x_{1})
\cdots{\cal Y}_{m}(w_{(\alpha_{m})},
x_{m})w_{(\beta_{m})}\rangle_{W^{\mu_{1}}}\mbar_{x^{n}_{i}=e^{n\log z_{i}},
i=1, \dots, m, n\in {\Bbb R}}$$ 
defined on the region $|z_{1}|>\cdots>|z_{m}|>0$
can be analytically extended to the {\it generalized rational
function}
$$\sum_{j=1}^{p}
g_{j}(z_{1}, \dots, z_{m})
\prod_{k=1}^{m}z_{k}^{-r^{(k)}_{j}}\prod_{k>l}(z_{1}-z_{2})^{-s^{(kl)}_{j}}.$$
\end{theo}

We omit the proof of this theorem  since 
it is completely analogous to the proof of
Theorem \ref{dual}.

For constructions and further study of intertwining operator algebras, see
\cite{H1}, \cite{H2}, \cite{H4} and \cite{H5}.

\renewcommand{\theequation}{\thesection.\arabic{equation}}
\renewcommand{\therema}{\thesection.\arabic{rema}}
\setcounter{equation}{0}
\setcounter{rema}{0}
\section{Genus-zero modular functors and 
algebras over them from intertwining operator algebras}

We construct in this section a genus-zero modular functor
from an intertwining operator algebra and prove that the intertwining 
operator algebra gives an algebra over the  partial operad underlying this
genus-zero modular functor, satisfying 
some additional conditions, including 
a certain generalized meromorphicity. 

Let $(W, {\cal A}, \{{\cal V}_{\alpha_{1}\alpha_{2}}^{\alpha_{3}}\},
{\bf 1}, \omega)$ be an intertwining operator algebra of central
charge $c$. We first construct holomorphic vector bundles ${\cal
M}^{W; \alpha_{0}}_{\alpha_{1}\cdots\alpha_{j}}(j)$ over $K(j)$ for
$j\in {\Bbb N}$ and $\alpha_{0}, \dots, \alpha_{j}\in {\cal A}$.  For
$j=0$ and $\alpha_{0}\in {\cal A}$, we take ${\cal M}^{W;
\alpha_{0}}(0)$ over $K(0)$ to be the trivial product line bundle
$K(0)\times {\Bbb C}$ over $K(0)$ if $\alpha_{0}=\epsilon$, and to be
the rank-$0$ bundle over $K(0)$ if $\alpha_{0}\ne \epsilon$.

In the case of $j=1$, for any $\alpha_{0}\in {\Bbb A}$,
we have a holomorphic multivalued operator-valued function
$e^{(\alpha_{0})}$ on $K(1)$ valued in 
$\mbox{\rm Hom}(W^{\alpha_{0}}, \overline{W}^{\alpha_{0}})$:
The values of $e^{(\alpha_{0})}$ at $(A^{(0)}, (a_{0}^{(1)}, A^{(1)}))\in
K(1)$ is
$$e^{(\alpha_{0})}_{(A^{(0)}, (a_{0}^{(1)}, A^{(1)}))}=
e^{-L^{-}_{A^{(0)}}}e^{-L^{+}_{A^{(1)}}}(a_{0}^{(1)})^{-L(0)}$$ 
where for a sequence $A=\{A_{j}\}_{j\in {\Bbb Z}_{+}}$
of complex numbers, 
\begin{eqnarray*}
L^{-}_{A}&=&\sum_{j\in {\Bbb Z}_{+}}A_{j}L(-j),\\
L^{+}_{A}&=&\sum_{j\in {\Bbb Z}_{+}}A_{j}L(j).
\end{eqnarray*}
This multivalued operator-valued function gives a holomorphic line bundle
${\cal
M}^{W; \alpha_{0}}_{\alpha_{0}}(1)$ over $K(1)$. 

Explicitly, this 
holomorphic bundle is constructed
as follows: Given any simply connected open subset of
$K(1)$, the restriction of the holomorphic 
multivalued operator-valued function above have single-valued branches and any 
two such branches are different from each other by a constant factor. 
These single-valued branches gives a line bundle over the open subset 
together with a holomorphic trivialization. 
Now consider  the intersection of two such open subsets. The connected
components of such an intersection must be simply connected. 
We consider one of such component. Through 
such an component, a single-valued
branch on one of the open subset discussed above has a unique analytic 
extension to a single-valued branch on the other open subset. 
Together with the trivializations of the line bundles over 
the two open subsets given by the the single-valued branches 
on these two open subsets, we obtain a transition function defined on this
component. Combining the transition functions on all the components of
the intersection, we obtain a transition function defined on 
the intersection. From the definition, we see that these transition functions
satisfy the required properties and thus give a holomorphic line
bundle ${\cal M}^{W; 
\alpha_{0}}_{\alpha_{0}}(1)$. 
For any $\alpha_{0}, \alpha_{1}\in {\cal A}$ such that $\alpha_{0}\ne 
\alpha_{1}$, we take the vector bundle ${\cal M}^{W; 
\alpha_{0}}_{\alpha_{1}}(1)$ to be the rank-$0$ bundle over $K(1)$.

We now construct ${\cal M}^{W; \alpha_{0}}_{\alpha_{1}\cdots\alpha_{j}}(j)$ 
for 
$j\ge 2$ and $\alpha_{0}, \dots, \alpha_{j}\in {\cal A}$.  
We need the notion of generalized-meromorphic function on $K(m)$, 
$m\in {\Bbb N}$. For $m\in {\Bbb N}$, 
a {\it generalized-meromorphic function on $K(m)$} is a multivalued map
from $K(m)$ to ${\Bbb C}$ which can be expressed as a finite sum  
whose summands are polynomials in $z_{1}, \dots, z_{m-1}$, $a_{0}^{(1)},
\dots, a_{0}^{(m)}$, $A^{(0)}_{j}, \dots, A^{(m)}_{j}$, $j\in {\Bbb Z}_{+}$,
divided by real powers of  $z_{1}, \dots, z_{m-1}$, 
$a_{0}^{(1)}, \dots, a_{0}^{(m)}$ and $z_{k}-z_{l}$, $k\ne l$. 
The points $z_{i}=0, \infty$,
$z_{k}=z_{l}$, $a_{0}^{(i)}=0, \infty$, $A^{(i)}_{j}=\infty$
are called the singularities of this generalized-meromorphic 
function on $K(m)$.
The real powers of $z_{k}-z_{l}$ appearing
in the generalized-meromorphic function
such that the corresponding term are not zero are
called the {\it orders of the singularity $z_{k}=z_{l}$}.
By definition, a meromorphic function on $K(m)$ defined in \cite{H3} is a
generalized-meromorphic function.

By the definition of intertwining operator
algebra, we know that 
for any
$\alpha_{0}, \dots, \alpha_{j}, \beta_{1}, \dots, \beta_{j-2}\in {\cal A}$,  
${\cal Y}_{1}\in {\cal 
V}_{\alpha_{1}\;\beta_{1}}^{\alpha_{0}}$, ${\cal Y}_{i}\in {\cal 
V}_{\alpha_{i}\;\beta_{i+1}}^{\beta_{i}}$, $i=2, \dots, j-2$, 
${\cal Y}_{j-1}\in {\cal 
V}_{\alpha_{j-1}\;\alpha_{j}}^{\beta_{j-2}}$, any $w_{(\alpha_{i})}\in 
W^{\alpha_{i}}$,
$i=1, \dots, j$, $w'_{(\alpha_{0})}\in (W^{\alpha_{0}})'$,  the series
\begin{eqnarray*}
\lefteqn{\langle e^{-L^{+\prime}_{A^{(0)}}}w'_{(\alpha_{0})}, 
{\cal Y}_{1}(e^{-L^{+}_{A^{(1)}}}(a_{0}^{(1)})^{-L(0)}
w_{(\alpha_{1})}, x_{1})\cdots{\cal Y}_{j-1}(e^{-L^{+}_{A^{(j-1)}}}
(a_{0}^{(j-1)})^{-L(0)}\cdot}\nno\\
&&\hspace{8em}\cdot w_{(\alpha_{j-1})},
x_{j-1})e^{-L^{+}_{A^{(j)}}}
(a_{0}^{(j)})^{-L(0)}w_{(\alpha_{j})}\rangle\mbar_{x_{i}=z_{i},
i=1, \dots, j-1}
\end{eqnarray*}
is absolutely convergent when $|z_{1}|>\cdots >|z_{j-1}|>0$. 
Using the associativity , commutativity and generalized rationality 
for intertwining operators (Theorems \ref{dual} and 
\ref{dual2}) and the weight-grading-restriction
condition, it is easy to see that
this multivalued analytic function defined 
on the region determined by $|z_{1}|>\cdots >|z_{j-1}|>0$ can be 
analytically extended
to a generalized-meromorphic function on $K(j)$. 
Thus for fixed 
$\alpha_{0}, \dots, \alpha_{j}, \beta_{1}, \dots, \beta_{j-2}\in {\cal A}$,  
${\cal Y}_{1}\in {\cal 
V}_{\alpha_{1}\;\beta_{1}}^{\alpha_{0}}$, ${\cal Y}_{i}\in {\cal 
V}_{\alpha_{i}\;\beta_{i+1}}^{\beta_{i}}$, $i=2, \dots, j-2$, 
${\cal Y}_{j-1}\in {\cal 
V}_{\alpha_{j-1}\;\alpha_{j}}^{\beta_{j-2}}$,  we obtain 
an multivalued operator-valued function $m_{{\cal Y}_{1}\otimes \cdots
\otimes {\cal Y}_{j-1}}$ on
$K(j)$ valued in $\hom(W^{\alpha_{1}}\otimes \cdots\otimes W^{\alpha_{j}}, 
\overline{W}_{\alpha_{0}})$ defined by
\begin{eqnarray*}
\lefteqn{m_{{\cal Y}_{1}\otimes \dots \otimes 
{\cal Y}_{j-1}}(Q(j, z, a, A))}\nno\\
&&= 
e^{-L^{-}_{A^{(0)}}}{\cal Y}_{1}(e^{-L^{+}_{A^{(1)}}}(a_{0}^{(1)})^{-L(0)}
\cdot, x_{1})\cdots\nno\\
&&\quad \cdots {\cal Y}_{j-1}
(e^{-L^{+}_{A^{(j-1)}}}(a_{0}^{(j-1)})^{-L(0)}\cdot,
x_{j-1})e^{-L^{+}_{A^{(j)}}}(a_{0}^{(j)})^{-L(0)}
\end{eqnarray*}
for
$$Q(j, z, a, A)=(z_{1}, \dots, z_{j-1}; A(z; 1), (a_{0}^{(1)}, {\bf 0}),
\dots, (a_{0}^{(j)}, {\bf 0}))\in K(j).$$
For fixed 
$\alpha_{0}, \dots, \alpha_{j}\in {\cal A}$, these 
multivalued operator-valued functions determine a holomorphic vector bundle 
${\cal M}_{\alpha_{1}\cdots \alpha_{j}}^{W; \alpha_{0}}(j)$
of finite rank as in the case of $j=1$ above.

Recall from \cite{H3} that for any $j\in {\Bbb N}$,
the $c/2$-power $\tilde{K}^{c}(j)$ of the determinant line bundle over $K(j)$
is the trivial product line bundle. Thus the tensor product of 
$\tilde{K}^{c}(j)$ and any holomorphic vector bundle over $K(j)$ is 
canonically isomorphic to the holomorphic vector bundle itself. 
Because of this
fact, we shall often identify 
${\cal M}_{\alpha_{1}\cdots \alpha_{j}}^{W; \alpha_{0}}(j)$ 
with
${\cal M}_{\alpha_{1}\cdots \alpha_{j}}^{W; \alpha_{0}}(j)\otimes 
\tilde{K}^{c}(j)$.

Taking direct sums of the holomorphic vector bundles we just constructed,
we obtain a sequence ${\cal M}^{W}$ 
of finite-rank holomorphic vector bundles 
$${\cal M}^{W}(j)=\oplus_{\alpha_{0}, \dots, \alpha_{j}\in {\cal A}}
{\cal M}_{\alpha_{1}\cdots \alpha_{j}}^{W; \alpha_{0}}(j),\;\;j\in 
{\Bbb N}.$$

Next we define the composition maps for ${\cal M}_{W}$.  As usual, we
use $(Q, F)$ to denote an element of ${\cal M}(j)$ for any $j\in {\Bbb
N}$ where $Q\in K(j)$ and $F$ is in the fiber over $Q$.  Let $(Q_{1},
F_{1})\in {\cal M}_{W}(j)$ and $(Q_{2}, F_{2})\in {\cal M}_{W}(k)$. Then
by definition, there are multivalued operator-valued functions which 
are linear
combinations of $m_{{\cal Y}^{(1)}_{1}, \dots, {\cal Y}^{(1)}_{j-1}}$
and $m_{{\cal Y}^{(2)}_{1}, \dots, {\cal Y}^{(2)}_{k-1}}$ such that
$F_{1}$ and $F_{2}$ are values  of these functions at $Q_{1}$ and
$Q_{2}$, respectively.  To define the composition maps, we need only
consider the case in which $F_{1}$ and $F_{2}$ are values of $m_{{\cal
Y}^{(1)}_{1}, \dots, {\cal Y}^{(1)}_{j-1}}$ and $m_{{\cal Y}^{(2)}_{1},
\dots, {\cal Y}^{(2)}_{k-1}}$ at $Q_{1}$ and $Q_{2}$, respectively.
Assume that for $1\le i\le j$, $Q_{1^{i}}\infty_{^{0}}Q_{2}$ exists.
Choose simply connected open subsets $U_{1}$ and $U_{2}$ of $K(j)$ and
$K(k)$ containing $Q_{1}$ and $Q_{2}$, respectively, such that the image
of the direct product of these open subsets under the sewing operation
is a simply connected open subset $U_{3}$ of $K(j+k-1)$ containing
$Q_{1^{i}}\infty_{^{0}}Q_{2}$. Using the associativity, commutativity
and skew-symmetry for intertwining operators instead of the
corresponding properties for vertex operators associated to a vertex
operator algebra, we can adapt the proof of Proposition 5.4.1 in
\cite{H3} to show that if we use  $(m_{{\cal Y}^{(1)}_{1}, \dots,
{\cal Y}^{(1)}_{j-1}})(Q_{1})$ and $(m_{{\cal Y}^{(2)}_{1},
\dots, {\cal Y}^{(2)}_{k-1}})(Q_{2})$ to denote the set of all values 
of  $m_{{\cal Y}^{(1)}_{1}, \dots,
{\cal Y}^{(1)}_{j-1}}$ and $m_{{\cal Y}^{(2)}_{1},
\dots, {\cal Y}^{(2)}_{k-1}}$ at $Q_{1}$ and $Q_{2}$,
respectively, then the set of
the contractions 
$$(m_{{\cal Y}^{(1)}_{1}, \dots,
{\cal Y}^{(1)}_{j-1}})(Q_{1})_{^{i}}*_{^{0}} (m_{{\cal Y}^{(2)}_{1},
\dots, {\cal Y}^{(2)}_{k-1}})(Q_{2})$$   
is in fact the set of values of a linear
combination of the multivalued operator-valued functions of the form $m_{{\cal
Y}_{1}, \dots, {\cal Y}_{j+k-2}}$ at $Q_{1^{i}}\infty_{^{0}}Q_{2}$. 
(Note that the central charge $c$ of the intertwining operator algebra $W$ 
enters in the coefficients of the linear combination.)

We
take a particular value of this linear combination as follows: Consider
the single-valued branches of the restrictions of $m_{{\cal
Y}^{(1)}_{1}, \dots, {\cal Y}^{(1)}_{j-1}}$ and $m_{{\cal Y}^{(2)}_{1},
\dots, {\cal Y}^{(2)}_{k-1}}$ to $U_{1}$ and $U_{2}$, respectively, such
that $F_{1}$ and $F_{2}$ are the restrictions of these branches at
$Q_{1}$ and $Q_{2}$, respectively. The contraction of these branches
gives a single-valued branch of the restriction of the linear
combination above to $U_{3}$. The value of this single-valued branch at
$Q_{1^{i}}\infty_{^{0}}Q_{2}$ is the particular value we want.  Denote
this particular value by $F_{3}$. We define $(Q_{1},
F_{1})_{^{i}}\infty^{W}_{^{0}}(Q_{2}, F_{2})$ to be
$(Q_{1^{i}}\infty_{^{0}}Q_{2}, F_{3})$. Using this operation, we can
define the substitution map $\gamma_{{\cal M}^{W}}$ for ${\cal M}^{W}$
similarly to the substitution map
for the vertex partial operads $\tilde{K}^{c}$, $c\in
{\Bbb C}$, in \cite{H3}. 

In ${\cal M}^{W}(1)$, we have an identity $I_{{\cal M}^{W}}$ which 
corresponds to the operator
whose restriction to $W^{\alpha}$ for 
$\alpha\in {\cal A}$ is
the identity operator $I_{W^{\alpha}}: W^{\alpha}\to W^{\alpha}$.
Also there is an obvious action of $S_{j}$ on ${\cal M}^{W}(j)$ for 
$j\in {\Bbb N}$.

We have:

\begin{theo}\label{w-mf}
The sequence ${\cal M}_{W}=\{{\cal M}_{W}(j)\}_{j\in {\Bbb N}}$ together
with the substitution $\gamma_{{\cal M}^{W}}$, the identity 
$I_{{\cal M}^{W}}$ and the actions of $S_{j}$, $j\in {\Bbb N}$,  is 
a rational genus-zero modular functor. 
\end{theo}
\pf
From the construction of ${\cal M}_{W}$ above, we see that the only
axiom need to be verified is the
operad-associativity. Let $(Q_{1}, F_{1})\in {\cal M}^{W}(j)$, 
$(Q_{2}, F_{2})\in {\cal M}^{W}(k)$ and 
$(Q_{3}, F_{3})\in {\cal M}^{W}(l)$. 
Assume that for $1\le i_{1}\le j$ and
$i_{1} \le i_{2}\le i_{1}+k-1$, both
$(Q_{1\;^{i_{1}}}\infty_{^{0}}Q_{2})_{^{i_{2}}}\infty_{^{0}}Q_{3}$
and $Q_{1\;^{i_{1}}}\infty_{^{0}}(Q_{2\;^{i_{2}-i_{1}-1}}\infty_{^{0}}Q_{3})$
exist. Choose simply connected open subsets $U_{1}\subset K(j)$,
$U_{2}\subset K(k)$ and $U_{3}\subset  K(l)$ such that 
for any elements $Q'_{1}\in U_{1}$, $Q'_{2}\in U_{2}$ and $Q'_{3}\in U_{3}$,
both $(Q'_{1\;^{i_{1}}}\infty_{^{0}}Q'_{2})_{^{i_{2}}}\infty_{^{0}}Q'_{3}$
and $Q'_{1\;^{i_{1}}}\infty_{^{0}}(Q'_{2\;^{i_{2}-i_{1}-1}}
\infty_{^{0}}Q'_{3})$
exist and such that the open subsets 
\begin{equation}\label{2.0}
(U_{1\;^{i_{1}}}\infty_{^{0}}U_{2})_{^{i_{2}}}\infty_{^{0}}U_{3}=
U_{1\;^{i_{1}}}\infty_{^{0}}(U_{2\;^{i_{2}-i_{1}-1}}\infty_{^{0}}U_{3})
\end{equation}
of $K(j+k+l-2)$ are simply connected.
Then using the definition of the composition for ${\cal M}^{W}$ and
commutativity and 
skew-symmetry for intertwining operators,
and adapting the arguments in the proof of
Proposition 5.4.1 in \cite{H3}, 
we can show that both 
\begin{equation}\label{2.1}
((Q_{1}, F_{1})_{^{i_{1}}}\infty^{W}_{^{0}}(Q_{2}, F_{2}))_{^{i_{2}}}
\infty^{W}_{^{0}}(Q_{3}, F_{3})
\end{equation}
and 
\begin{equation}\label{2.2}
(Q_{1}, F_{1})_{^{i_{1}}}
\infty^{W}_{^{0}}((Q_{2}, F_{2})_{^{i_{2}-i_{1}+1}}
\infty^{W}_{^{0}}(Q_{3}, F_{3})
\end{equation}
are restrictions at $(Q_{1\;^{i_{1}}}
\infty_{^{0}}Q_{2})_{^{i_{2}}}\infty_{^{0}}Q_{3}$
and $Q_{1\;^{i_{1}}}\infty_{^{0}}(Q_{2\;^{i_{2}-i_{1}-1}}\infty_{^{0}}Q_{3})$,
respectively,
of a single-valued branch on (\ref{2.0}) of 
a linear combination of the multivalued operator-valued functions of the form 
$m_{{\cal Y}_{1}, \dots, {\cal Y}_{j+k+l-3}}$. Since 
$$(Q_{1\;^{i_{1}}}
\infty_{^{0}}Q_{2})_{^{i_{2}}}\infty_{^{0}}Q_{3}=
Q_{1\;^{i_{1}}}\infty_{^{0}}(Q_{2\;^{i_{2}-i_{1}-1}}\infty_{^{0}}Q_{3}),$$
(\ref{2.1}) and (\ref{2.2}) are equal, proving the operad-associativity.
\epfv

In fact the constructions and proof above also give 
an algebra over the 
$({\cal M}^{W})^{\times}(1)|_{{\Bbb C}^{\times}}$-rescalable
partial operad ${\cal M}^{W}$ 
satisfying
certain additional properties
(recall the definition of algebra over 
a rescalable partial operad in \cite{HL1}, \cite{HL2} and \cite{H3}). 
For simplicity, by an {\it algebra over a rational genus-zero modular 
functor ${\cal M}$}, we shall mean an algebra
over the underlying 
${\cal M}^{\times}(1)|_{{\Bbb C}^{\times}}$-rescalable
partial operad of ${\cal M}$. To describe such an algebra, we first need 
discuss a certain class of
irreducible modules for ${\cal M}^{\times}(1)|_{{\Bbb C}^{\times}}$.

Since ${\cal M}(1)|_{{\Bbb C}^{\times}}$ is in fact 
a line bundle over ${\Bbb C}^{\times}$,  the restriction 
${\cal M}^{\times}(1)|_{{\Bbb R}_{+}}$ of
${\cal M}^{\times}(1)|_{{\Bbb C}^{\times}}$ to the positive real numbers
can be identified, after choosing a suitable section,
with the space of pairs $(\alpha, \mu)$,
$\alpha\in {\Bbb R}_{+}$, $\mu\in {\Bbb C}^{\times}$. Note that 
${\cal M}^{\times}(1)|_{{\Bbb R}_{+}}$ is a subgroup of 
${\cal M}^{\times}(1)|_{{\Bbb C}^{\times}}$ and the product is given by
$$\gamma_{\cal M}((\alpha_{1}, \mu_{1}); (\alpha_{2}, \mu_{2}))
=(\alpha_{1}\alpha_{2}, \mu_{1}\mu_{2})$$
for $\alpha_{1}, \alpha_{2}\in {\cal R}_{+}$ and
$\mu_{1}, \mu_{2}\in {\Bbb C}^{\times}$. In particular, the subset of 
${\cal M}^{\times}(1)|_{{\Bbb R}_{+}}$ consisting of pairs $(\alpha, 1)$,
$\alpha\in {\Bbb R}_{+}$, is a subgroup of ${\cal M}^{\times}(1)
|_{{\Bbb R}_{+}}$
isomorphic to the group ${\Bbb R}_{+}$. Let $W$ be an irreducible module
for ${\Bbb R}_{+}$. Then $W$ must be one-dimensional and
there is a complex number $n$ such that the action of ${\Bbb R}_{+}$
on $W$ is given by 
\begin{eqnarray*}
{\Bbb R}_{+}&\to &\mbox{\rm End}\; W\nno\\
\alpha&\mapsto &\alpha^{n}.
\end{eqnarray*}

We shall call 
$n$ the {\it weight} of $W$. Given any irreducible 
module $W$ for ${\Bbb R}_{+}$ of weight $n\in {\Bbb C}$ and any flat
holomorphic connection on ${\cal M}(1)|_{{\Bbb C}^{\times}}$,  we define an
irreducible 
${\cal M}(1)|_{{\Bbb C}^{\times}}$-module structure on $W$ as follows:
Using  the 
flat holomorphic connection on ${\cal M}(1)|_{{\Bbb C}^{\times}}$,
we identify ${\cal M}(1)|_{{\Bbb C}^{\times}}$ with the set of
pairs of the form $(l, \mu)$ where $l\in {\Bbb C}$ satisfying 
$0\le \Im{(l)}< 2\pi$, and $\mu \in {\Bbb C}$. Then the 
${\cal M}(1)|_{{\Bbb C}^{\times}}$-module structure is defined by
\begin{eqnarray*}
{\cal M}(1)|_{{\Bbb C}^{\times}}&\to &\mbox{\rm End}\; W\nno\\
(l, \mu)&\mapsto &\mu e^{nl}.
\end{eqnarray*}
We call such an irreducible ${\cal M}(1)|_{{\Bbb C}^{\times}}$-module
an {\it irreducible ${\cal M}(1)|_{{\Bbb C}^{\times}}$-module induced from
an irreducible ${\Bbb R}_{+}$-module of  weight $n$}.

Recall that the underlying vector space of 
an algebra over an ${\cal M}(1)|_{{\Bbb C}^{\times}}$-rescalable
partial operad is a direct sum of irreducible 
${\cal M}^{\times}(1)|_{{\Bbb C}^{\times}}$-modules (see \cite{HL1},
\cite{HL2} and Appendix C in \cite{H3}). 
We have the following notion:

\begin{defi}\label{g-m-a}
{\rm An {\it
generalized-meromorphic algebra  over a genus-zero modular 
functor ${\cal M}$} is an algebra $(W, V, \Upsilon)$ over ${\cal M}$
satisfying the following conditions:

\begin{enumerate}

\item $W=\coprod_{n\in {\Bbb R}}W_{(n)}$ where for any 
$n\in {\Bbb R}$, $W_{(n)}$ is an
irreducible 
${\cal M}(1)|_{{\Bbb C}^{\times}}$-module induced from an
irreducible ${\Bbb R}_{+}$-module of weight $n$ (note that we have a
flat holomorphic connection on ${\cal M}(1)|_{{\Bbb C}^{\times}}$).

\item  $W_{(n)}=0$ for $n$ sufficiently small.

\item For any $j\in {\Bbb N}$, $\Upsilon_{j}$ is linear on any fiber
of ${\cal M}(j)$.

\item \label{mero} Let $j\in {\Bbb N}$ and let $U$ be any simply connected open
subset of $K(j)$. Then for any flat section $\psi: U\to 
{\cal M}(j)|_{U}$ of the restriction ${\cal M}(j)|_{U}$ of
the vector bundle ${\cal M}(j)$ and any $w'\in W'$,
$w_{1}, \dots, w_{j}\in W$, the function
$$Q\to \langle w', (\Upsilon_{j}(\psi(Q))) 
(w_{1}\otimes \cdots
\otimes w_{j})\rangle$$ 
on $U$ can be analytically extended to a generalized-meromorphic function
on $K(j)$ satisfying the following property:
If $z_{k}$ and $z_{l}$ are the $k$-th and
$l$-th punctures of $Q\in K(j)$ respectively, 
then for any $w_{k}$ and $w_{l}$
in $W$ there
exists a positive integer $N(w_{k}, w_{l})$
 such that for any $w'\in W'$,
$w_{i}\in W$,
$i\ne k, l$, the  
orders of the singularity $z_{k}=z_{l}$ (we use the convention
$z_{j}=0$) of $$\langle w', (\Upsilon_{j}(\psi(Q))) 
(w_{1}\otimes \cdots
\otimes w_{j})\rangle$$
is  less than $N(w_{k}, w_{l})$.

\end{enumerate}}
\end{defi}

We shall denote the generalized-meromorphic algebra over ${\cal M}$
above by $(W, V, \Upsilon)$.

Note that from the generalized-meromorphic functions in Axiom
\ref{mero} above, we can construct holomorphic vector bundles using the
method in the construction of a modular functor from an intertwining
operator algebra above. Clearly these holomorphic vector bundles are
isomorphic to ${\cal M}(j)$, $j\in {\Bbb N}$. But in general they are
certainly not the same as ${\cal M}(j)$, $j\in {\Bbb N}$. To formulate
the isomorphism theorem (Theorem \ref{main} 
in the next section precisely, we need the
following notion: 

\begin{defi}
{\rm A {\it 
canonically-generalized-meromorphic algebra} 
over a rational genus-zero modular functor ${\cal M}$ is a generalized
meromorphic algebra  $(W, V, \Upsilon)$ over ${\cal M}$ such that 
 the vector bundles on ${\cal M}(j)$, $j\in
{\Bbb N}$, is  the same (not just ispomorphic to) as 
the ones constructed from the generalized-meromorphic functions
given by Axiom \ref{mero} in Definition \ref{g-m-a}.}
\end{defi}

{\it Homomorphisms} and {\it isomorphisms} between
generalized-meromorphic algebras over genus-zero modular functors are
defined in the obvious way.

\begin{rema}\label{weakly}
{\rm If  the fourth axiom in the above definition is replaced
by the holomorphicity of $\Upsilon$, then we do not need the flat
holomorphic connections $\nabla(j)$ 
on ${\cal M}(j)$, $j\in {\Bbb N}$. An algebra $(W, V, \Upsilon)$
over a genus-zero
modular functor ${\cal M}$
satisfying the first three axioms in the definition
above and the holomorphicity of $\Upsilon$  is called a
genus-zero weakly holomorphic conformal field theory over ${\cal M}$
(see \cite{H2}, \cite{S1} and \cite{S2}). 
The fourth axiom above clearly implies the holomorphicity of
$\Upsilon$. Thus the underlying algebra over ${\cal M}$ of 
a generalized-meromorphic algebra over ${\cal M}$ is
a genus-zero weakly holomorphic conformal field theory over ${\cal
M}$.}
\end{rema}

Let $(W, {\cal A}, \{{\cal V}_{\alpha_{1}\alpha_{2}}^{\alpha_{3}}\}, 
{\bf 1}, \omega)$ be an intertwining operator algebra. Let $V$ be the 
vertex operator subalgebra of $W^{e}$ generated by ${\bf 1}$ and $\omega$.
By Theorem 
\ref{w-mf}, we have a rational genus-zero modular functor ${\cal M}^{W}$.
We define a morphism of partial operad
$$\Upsilon^{W}: {\cal M}^{W}\to {\cal H}_{W, V}$$
as follows:
Given any element $(Q, F)$ of ${\cal M}^{W}_{j}$ for $j\in {\Bbb N}$,
by definition, $F$ is a value of a linear combination of
generalized-meromorphic functions
on $K(j)$ of the form $m_{{\cal Y}_{1}, \dots, {\cal Y}_{j-1}}$. 
So $F$ is an element of ${\cal H}_{W, V}(j)$. We define 
$\Upsilon^{W}_{j}((Q, F))=F$. 

We have:

\begin{theo}\label{->}
The quadruple  $(W, V, \Upsilon^{W})$ 
is a generalized-meromorphic algebra over the genus-zero 
modular functor ${\cal M}^{W}$. In particular, $(W, V, \Upsilon^{W})$
is a genus-zero weakly holomorphic conformal field theory over ${\cal M}^{W}$.
\end{theo}
\pf
Using the duality properties for intertwining operator algebras
instead of the corresponding
properties for vertex operator algebras, we can adapt 
the argument in the proof of Proposition 5.4.1 of \cite{H3} to
prove that $(W, V, \Upsilon^{W})$ is an algebra over ${\cal M}^{W}$. 
The generalized meromorphicity is clear from the definition of 
$\Upsilon$ and the constructions of the holomorphic flat connection on
${\cal M}^{W}(j)$, $j\in {\Bbb N}$, in Section 1 and this section. 
The second conclusion follows from Remark \ref{weakly}.
\epfv

Combining Theorems \ref{->} and \ref{braid}, we obtain:

\begin{corol}
Let $(W, {\cal A}, \{{\cal V}_{\alpha_{1}\alpha_{2}}^{\alpha_{3}}\}, 
{\bf 1}, \omega)$ be an intertwining operator algebra and
${\cal V}_{W}$ the fiber over any point in $F(j)\subset K(j)$ of the vector
bundle ${\cal M}_{W}(j)$. Then ${\cal V}$ has a natural structure of a
representation of the framed braid group  on $j$ strings. In particular, 
${\cal V}$ has a natural structure of a
representation of the braid group $B_{j}$ on $j$ strings.
\end{corol}

\renewcommand{\theequation}{\thesection.\arabic{equation}}
\renewcommand{\therema}{\thesection.\arabic{rema}}
\setcounter{equation}{0}
\setcounter{rema}{0}
\section{The equivalence theorem}

In this section, we first show that a generalized-meromorphic algebra
over a rational genus-zero modular functor gives an intertwining operator 
algebra. The main result of the present paper follows easily.

Let $(W, V, \Upsilon)$ be a generalized-meromorphic algebra
over a rational genus-zero modular functor ${\cal M}$. By the definition of 
rational genus-zero modular functor, we already have a finite set ${\cal A}$.
Let $\epsilon$ be the element of ${\cal A}$ 
such that the identity $I_{\cal M}$ is
in ${\cal M}^{\epsilon}_{\epsilon}(1)$. 
For any $a_{1}, a_{2}, a_{3}\in {\cal A}$, 
by the definition of generalized meromorphicity, 
the image under $\Upsilon_{2}$
of any flat section of ${\cal M}_{a_{1}a_{2}}^{a_{3}}$ over a 
simply connected open subset of ${\cal M}(2)$ is an element of
${\cal H}_{W, V}(2)$ such that its matrix elements can be analytically
extended to generalized-meromorphic functions. Considering the 
subset $M^{1}$ of $K(2)$, we see that the image of such a 
flat section in particular gives an multivalued operator-valued function 
on $M^{1}$, which can be identified with ${\Bbb C}^{\times}$. Such an 
multivalued operator-valued function on ${\Bbb C}^{\times}$ gives a
map  from $W^{a_{1}}\otimes W^{a_{2}}\to W^{a_{3}}\{x\}$.
We 
define ${\cal V}_{a_{1}a_{2}}^{a_{3}}$ to be the space of 
all maps obtained in this way. 
The vacuum ${\bf 1}$ and the Virasoro element $\omega$
are defined by
\begin{equation}
{\bf 1}=\Upsilon_{0}(({\bf 0}; 1));
\end{equation}
and
$$
\omega=-{\displaystyle
\frac{d}{d\varepsilon}}\nu((A(\varepsilon; 2); 1))\lbar_{\varepsilon =0}.
$$

We have:

\begin{theo}\label{inverse}
Let $(W, V, \Upsilon)$ be a generalized-meromorphic algebra
over a rational genus-zero modular functor ${\cal M}$
of central charge $c$. Then 
$(W, {\cal A}, {\cal V}_{a_{1}a_{2}}^{a_{3}}, {\bf 1}, \omega)$ is
an intertwining operator algebra of central charge $c$.
\end{theo}
\pf 
The proof of this result is completely analogous to the 
proof of Proposition 5.4.4 in \cite{H3} except for the commutative associative
structure on the vector space $A$ spanned by elements of ${\cal
A}$. This structure can be constructed as follows: For any
$\alpha_{1}, \alpha_{2}, \alpha_{3}\in {\cal A}$, let ${\cal
N}_{\alpha_{1}\alpha_{2}}^{\alpha_{3}}$ be the rank of 
${\cal M}_{\alpha_{1}\alpha_{2}}^{\alpha_{3}}$. We define the product
in $A$ by 
$$\alpha_{1}\cdot \alpha_{2}=\sum_{\alpha_{3}\in {\cal A}}{\cal
N}_{\alpha_{1}\alpha_{2}}^{\alpha_{3}}\alpha_{3}.$$
Then clearly this product is commutative and associative. 
From the axioms for rational genus-zero modular functors, it is easy
to see that ${\cal
N}_{\epsilon\alpha_{2}}^{\alpha_{3}}$ is $1$ when
$\alpha_{2}=\alpha_{3}$ and is $0$ when $\alpha_{2}\ne \alpha_{3}$. So 
$A$ has an identity  $\epsilon$. 
\epfv

Combining the constructions and results in the preceding and the present
sections, we obtain the following main result of this paper:

\begin{theo}\label{main} 
The category of intertwining operator algebras of central
charge $c$ and the category of canonically-generalized-meromorphic algebras 
 over rational genus-zero modular functors of central
charge $c$
are isomorphic.
\end{theo}
\pf
In the preceding, we have already constructed 
a functor from the category of  intertwining operator algebras of central
charge $c$ to the category of 
canonically-generalized-meromorphic algebras 
 over rational genus-zero modular functors of central
charge $c$. Theorem \ref{inverse} above gives
a functor from the category of canonically-generalized-meromorphic algebras 
over rational genus-zero modular functors of central
charge $c$
to the category of  intertwining operator algebra of central
charge $c$. The argument in the proof of the equivalence 
theorem for vertex operator algebras in Section 5.4 of \cite{H3}
can be adapted in the obvious way to show that the compositions of
these two functors are  the identity functors.
\epfv

{\small \sc Laboratory of Mathematics for Nonlinear Problems, 
Fudan University, Shanghai, 
China}

{\small and}

{\small \sc Department of Mathematics, Rutgers University,
110 Frelinghuysen Rd., Piscataway, NJ 08854-8019}

{\em E-mail address}: yzhuang@math.rutgers.edu


\begin{thebibliography}{FLM2}


\bibitem[B]{B}
R.~E.~Borcherds,
Vertex algebras, Kac-Moody algebras, and the Monster,
{\em Proc. Natl. Acad. Sci. USA} {\bf 83} (1986), 3068--3071.

\bibitem[FLM]{FLM}
I. B. Frenkel, J. Lepowsky and A. Meurman, A natural representation of
the Fischer-Griess monster with the modular function $J$ as character,
{\em Proc. Natl. Acad. Sci. USA} {\bf 81} (1984), 3256--3260.

\bibitem[H1]{H-1}
Y.-Z. Huang,
{\em On the geometric interpretation of vertex operator algebras},
Ph.D thesis, Rutgers University, 1990.

\bibitem[H2]{H0}
Y.-Z. Huang,
Geometric interpretation of vertex operator algebras,
{\em Proc. Natl. Acad. Sci. USA} {\bf 88} (1991), 9964--9968.

\bibitem[H3]{H1}
Y.-Z. Huang, A theory of tensor products for module categories
for a vertex operator algebra, IV, {\em J. Pure Appl. Alg.} {\bf 100}
(1995), 173-216.


\bibitem[H4]{H1.5}
Y.-Z. Huang, A nonmeromorphic extension of the moonshine module 
vertex operator algebra, in: {\em Moonshine, the Monster
and related topics, Proc. Joint Summer Research Conference, 
Mount Holyoke, 1994,} ed. C. Dong and G. Mason, 
Contemporary Math., Vol. 193, Amer. Math. Soc., Providence, 1996, 
123--148.


\bibitem[H5]{H2}
Y.-Z. Huang, Intertwining operator algebras, genus-zero modular
functors and genus-zero conformal field theories,  in: {\em Operads:
Proceedings of Renaissance Conferences}, ed. J.-L. Loday,
J. Stasheff, and A. A. Voronov, Contemporary Math., Vol. 202,
Amer. Math. Soc., Providence, 1997,  335--355.


\bibitem[H6]{H3}
Y.-Z. Huang, {\em Two-dimensional conformal geometry and vertex operator
algebras}, Progress in Mathematics, Vol. 148,
Birkh\"{a}user, Boston, 1997.

\bibitem[H7]{H4}
Y.-Z. Huang, A Jacobi identity for intertwining operator algebras,
to appear.


\bibitem[H8]{H5}
Y.-Z. Huang, Intertwining operator algebras, in preparation.


\bibitem[HL1]{HL1}
Y.-Z. Huang and J. Lepowsky, Vertex operator algebras and operads,
{\em The Gelfand Mathematical Seminar, 1990--1992}, ed. L. Corwin, I.
Gelfand and J. Lepowsky, Birkh\"{a}user Boston, 1993, 145--161.

\bibitem[HL2]{HL2}
Y.-Z.~Huang and J.~Lepowsky,
\newblock  Operadic formulation of the notion of vertex operator algebra,
in: {\em Mathematical Aspects of Conformal and Topological Field
Theories and Quantum Groups, Proc. Joint Summer Research Conference,
Mount Holyoke, 1992}, ed. P. Sally, M. Flato, J.  Lepowsky, N.
Reshetikhin and G. Zuckerman,
\newblock Contemporary Math., Vol. 175,
 Amer. Math. Soc., Providence, 1994, 131-148.

\bibitem[HL3]{HL3}
Y.-Z. Huang and J. Lepowsky, A theory of tensor products for module
categories for a vertex operator algebra, I, {\em Selecta
Mathematica, New Series} {\bf 1} (1995), 699-756.

\bibitem[HL4]{HL4}
Y.-Z. Huang and J. Lepowsky, A theory of tensor products for module
categories for a vertex operator algebra, II, {\em Selecta
Mathematica, New Series} {\bf 1} (1995), 757--786.

\bibitem[HL5]{HL5}
Y.-Z. Huang and J. Lepowsky, Tensor products of modules for a vertex
operator algebras and vertex tensor categories, in:
     {\em Lie Theory and Geometry,
in honor of Bertram Kostant,}
ed. R. Brylinski, J.-L. Brylinski, V. Guillemin, V. Kac,
Birkh\"{a}user Boston, 1994, 349--383.

\bibitem[HL6]{HL6}
Y.-Z. Huang and J. Lepowsky, A theory of tensor products for module
categories for a vertex operator algebra, III, {\em J. Pure
Appl. Alg.} {\bf 100} (1995),  141-171.





\bibitem[S1]{S1}
G. B.~Segal, Two-dimensional conformal field theories and modular
functors, in: {\em Proceedings of the IXth International Congress on
Mathematical Physics, Swansea, 1988},
Hilger, Bristol, 1989, 22--37.

\bibitem[S2]{S2}
G. B.~Segal,
The definition of conformal field theory,
preprint.


\end{thebibliography}
\end{document}